\documentclass[11pt]{article}
\usepackage{amsfonts}
\usepackage{amsfonts,amssymb,amsmath,amsbsy,amsthm,amscd}
\usepackage[utf8]{inputenc} %my
\usepackage[english]{babel}   % load Babel setup for English
\usepackage{indentfirst} %my
\begin{document}
%======================================================================
%\cleanbegin
%======================================================================

\begin{center}
ASYMPTOTIC RELIABILITY ANALYSIS OF THE SYSTEM 
WITH $N$ WORKING ELEMENTS AND REPAIRING DEVICE
\end{center}

\begin{center}
Golovastova E. A.
\end{center}

\textbf{Annotation.} This paper deals with system with $n$ identical elements and one repairing device. While one element working other ones stay in reserve. The distribution of element working and repairing times are supposed to be exponential. Here we obtain the asymptotic distribution of the system lifetime under conditions of its high reliability.


%------------------------------------------------------------------

\section{Introduction}
We consider the following model of the mathematical reliability theory [1]. The system consists of $ n $ identical elements, one of which works, while the others stay in reserve. Also we have one repairing device
The element working time $\eta$ and its recovery time $\xi$ have an exponential distribution. The device can restore only one element at an each moment of time, and if there are another broken elements, they stay in a queue. The system crashes at the time $\tau$, when all $n$ elements brakes down.

The most important characteristics of this  model and its  can be obtained by relying on the results from [1], [2]. Since real technical systems in our days are sufficiently reliable, there arises the problem  of an asymptotic analysis of systems with high reliability.

This paper established that under the condition of element fast repairing: $ E \eta \to 0 $, then the distribution of the normalized random variable $ \ ( \frac {E \xi} {E \eta})^ {n-1} \tau $ converges to exponential with parameter $E \xi$.

\section{Model description}
We consider a system consisting of $ n $ identical elements, one of which  performs certain functions, and the rest ones are in reserve. Working element can brake down and then it is replaced by one of operable. Repairing device is reliable. The system crashes at the moment when all $ n $ elements become inoperative. We will find the distribution of system operating time and its asymptotic behaviour when the elements are quickly restored or highly reliable, or their working time is long.

We assume that all the random variables that determine system functioning are independent and have an exponential distribution. Namely, elements working times are exponentially distributed with the parameter $ \lambda $, and the recovery times of the broken elements  are exponentially distributed with the parameter $ \mu $. The system state is determine by the number of broken elements $ q (t) $ at the time $ t $. When $ 0 <q (t) <n $, then one of the elements is working, one is restored by device, and $ q (t) -1 $ of the broken elements are waiting for recovery. Under such assumptions $ q (t) $ is the Markov chain with the absorbing state $ \{n \} $ and the absorption time:
$$
\tau = inf \{ t>0: \ q(t)=n \},
\eqno(1)
$$
and we suppose, that $q(0)=0$.

We will show, that:
$$
\lim_{\mu \rightarrow \infty} E \ e^{-s(\frac{\lambda}{\mu})^{n-1}\tau} = \frac{1}{\lambda+s}.
\eqno(2)
$$

By definition, $f(x)\sim g(x)$, if $\lim_{x \to 0} \frac{f(x)}{g(x)} = 1$.

\section{System operating time distribution}
Define:
$$
P_{j}(t)= P(q(t) = j, \ \tau > t | q(0) = 0 ),
$$
The following Kolmogorov system of differential equations is valid:
$$
P_{0}'(t) = -\lambda P_{0}(t) + \mu P_{1}(t), 
$$
$$
P_{j}'(t) = -(\lambda + \mu) P_{j}(t)+ \lambda P_{j-1}(t) + \mu P_{j+1}(t),\quad 0<j<n-1 
\eqno(3) 
$$
$$
P_{n-1}'(t) = -(\lambda + \mu) P_{n-1}(t) + \lambda P_{n-2}(t)
$$
with initial conditions: $P_{0}(0)=1$, $P_{j}(0)=0$, $j\neq 0$.

$\varphi_{j} (s) = \int_0^\infty e^{-st}P_{j}(t) \ dt$, $Re \ s \geq 0$ ~---~  Laplace transforms of system state probabilities.

We notice, that for small $h$:
$$
P(t < \tau \leq t+h)= \lambda h P_{n-1}(t) + \bar{\bar{o}}(h),
$$
so,
$$
\frac{d \ P(\tau \leq t)}{dt}=\lambda \ P_{n-1}(t).
\eqno(4) 
$$
By (4) we get, that $\lambda P_{n-1}(t)$ is the distribution density of $\tau$. That's why the convergence:
$$
P( (\lambda / \mu)^{n-1}\tau > t) \to e^{-\lambda t} \ \text{when} \ \mu \to \infty,
$$
is equivalent to convergence:
$$
\varphi_{n-1}((\lambda / \mu)^{n-1}s) \to \frac{1}{\lambda +s} \ \text{when} \ \mu \to \infty.
$$

\section{Asymptotic analysis of the operating time assuming its high reliability}
For $2$ elements in the system, we have the following equations for Laplace transforms of system state probabilities:
$$
s \ \varphi_0 - 1 = - \lambda \varphi_0 + \mu \varphi_1
$$
$$
s \ \varphi_1 = - (\lambda +\mu)\varphi_1 + \lambda \varphi_0.
$$
Due to brevity, we will sometimes omit the dependence on the argument below.
$$
\varphi_1(s) = \frac{\lambda}{s^2 + s(2\lambda + \mu) + \lambda^2}.
$$
Denote $\varepsilon = \frac{\lambda}{\mu}$; $\varepsilon \to 0$, when $\mu \to \infty$. So,
$$
\varphi_1(\varepsilon s) = \frac{\lambda}{\varepsilon^2 s^2 + \frac{\lambda}{\mu}s(2\lambda + \mu) + \lambda^2} \sim \frac{1}{\lambda + s}.
$$
For $3$ elements in the system, the correspondent equations are:
$$
s \ \varphi_0 - 1 = - \lambda \varphi_0 + \mu \varphi_1
$$
$$
s \ \varphi_1 = - (\lambda +\mu)\varphi_1 +\mu \varphi_2 + \lambda \varphi_0
$$
$$
s \ \varphi_2 = - (\lambda +\mu)\varphi_2 + \lambda \varphi_1,
$$
and the desired Laplace transform function is:
$$
\varphi_2(s) = \frac{\lambda^2}{(s+\lambda)^3 + 2\mu s^2 + 2 \mu\lambda s+\mu^2s}.
$$
So, when $\varepsilon \to 0$:
$$
\varphi_2(\varepsilon^2 s) \sim \frac{\lambda^2}{\lambda^3+\lambda^2s}=\frac{1}{\lambda+s}.
$$

Further, we assume that $n>3$.

For an arbitrary elements number $n$ in the system the equations for Laplace transform functions of the system state probabilities has the following form:
$$
\varphi_1 = \frac{\lambda+s}{\mu}\varphi_0-\frac{1}{\mu}
\eqno(5) 
$$
$$
\varphi_{j+1} = \frac{\lambda+\mu+s}{\mu}\varphi_j-\frac{\lambda}{\mu}\varphi_{j-1} \quad 0<j<n-1
\eqno(6) 
$$
$$
\varphi_{n-1} = \frac{\lambda}{\lambda+\mu+s}\varphi_{n-2}.
\eqno(7) 
$$

Expression (6) is the linear recurrence sequence with border conditions (5) and (7)[3]. So, according to the formula for the general term of such sequence, we have:
$$
\varphi_{j}(s)=A(s)q_1(s)^j + B(s)q_2(s)^j,
$$
where $q_1(s)$ and $q_2(s)$ are the roots of characteristic polynomial:
$$
q(s)^2 - \frac{\lambda+\mu+s}{\mu} q(s) + \frac{\lambda}{\mu}=0.
$$
$$
q_{1,2}(s)= \frac{(\lambda+\mu+s) \pm \sqrt{(\lambda+\mu+s)^2-4\lambda\mu} }{2\mu}.
$$
Here we notice, that:
$$
q_2 - q_1 = -\frac{\sqrt{(\lambda+\mu+s)^2-4\lambda\mu} }{\mu}, \quad q_1q_2=\frac{\lambda}{\mu}.
\eqno(8)
$$
$$
q_1(s) \sim 1+ \frac{s}{\mu}, \quad
q_2(s) \sim \frac{\lambda}{\mu}\ \text{when} \ s \to 0.
\eqno(9)
$$
From (5) and (7), we can get:
$$
A(s)= \frac{q_2^{n-2}(q_2(\lambda+\mu+s)-\lambda)}{q_2^{n-2}(\lambda-q_2(\lambda+\mu+s))(\mu q_1-(\lambda+s))+q_1^{n-2}(\mu q_2-(\lambda+s))(q_1(\lambda+\mu+s)-\lambda)}
$$
$$
B(s)=\frac{q_1^{n-2}(\lambda - q_1 (\lambda+\mu+s))}{q_2^{n-2}(\lambda-q_2(\lambda+\mu+s))(\mu q_1-(\lambda+s))+q_1^{n-2}(\mu q_2-(\lambda+s))(q_1(\lambda+\mu+s)-\lambda)}
$$
Using (8), we get:
$$
\varphi_{n-1}(s)= \frac{-\lambda\frac{\lambda}{\mu}^{n-2}\sqrt{\frac{s^2}{\mu^2}+2s(\frac{\lambda}{\mu^2} + \frac{1}{\mu}) +\frac{\lambda^2}{\mu^2} +1}}{q_2^{n-2}(\lambda-q_2(\lambda+\mu+s))(\mu q_1-(\lambda+s))+q_1^{n-2}(\mu q_2-(\lambda+s))(q_1(\lambda+\mu+s)-\lambda)}
$$
Denote $\varphi_{n-1}(s)= \frac{R(s)}{V(s)}$, and then, using (9), we find the asymptotic of $\varphi_{n-1}(\varepsilon^{n-1} s)$ when $\mu \to \infty$:
$$
R(\varepsilon^{n-1} s)\sim -\lambda\frac{\lambda}{\mu}^{n-2},
$$
$$
V(\varepsilon^{n-1} s) \sim \frac{\lambda}{\mu}^{n-2}(-\frac{\lambda^2}{\mu}+\frac{\lambda}{\mu}\varepsilon^{n-1} s)(\mu - \lambda - \varepsilon^{n-1} s) + (\mu + \varepsilon^{n-1} s)(-\varepsilon^{n-1} s)
$$
$$
\sim -\frac{\lambda}{\mu}^{n-2}\lambda^2 - \mu \varepsilon^{n-1} s = - \frac{\lambda}{\mu}^{n-2} (\lambda^2 + \lambda s).
$$
And so then again, when $\varepsilon \to 0$:
$$
\varphi_{n-1}(\varepsilon^{n-1} s) \sim \frac{1}{\lambda+s}.
$$
\\
\\
\begin{center}
\textbf{Literature}
\end{center}

1. Gnedenko B.V., Belyaev Yu.K., Soloviev A.D.
Mathematical methods in the theory of reliability. M.: Nauka, 1965.

2. Hutson V.C.L, Pym J.S.
Applications of Functional Analysis and Operator Theory. Academic Press, London, New York, Toronto, 1980, 432p. %{\bf22}.

3. Markushevich A.I. Recurrence sequences. Fizmatlit, 1950.

\end{document}